\def\version{0.70}

\def\journal{JA}
\def\titlep{$R$-matrices and the Yang-Baxter equation 
on GNS representations of C$^{*}$-bialgebras
\footnote{This is the revision of the previous version.}
}
\documentclass[11pt]{article}
\usepackage{graphicx,ifthen}
\usepackage{amssymb}

\font\germ=eufm10 at12pt

\def\goth#1{\hbox{\germ#1}}

\newcommand{\qed}{\hbox{\rule[-2pt]{3pt}{6pt}}}
\newcommand{\qedh}{\hfill\qed \\}

\newcommand{\vv}{\vspace{.3in}}

\def\labelenumi{\theenumi}
\def\theenumi{\arabic{enumi}}
\def\labelenumi{\theenumi}
\def\theenumi{\Alph{enumi}}
\renewcommand{\theenumi}{\Alph{enumi}}
\def\labelenumi{\theenumi}
\def\theenumi{\arabic{enumi}}
\def\labelenumi{\theenumi}
\def\theenumi{{\rm (\roman{enumi})}}

\newtheorem{Thm}{Theorem}[section]

\newtheorem{rem}[Thm]{Remark}

\newtheorem{ex}[Thm]{Example}
\newtheorem{defi}[Thm]{Definition}
\newtheorem{lem}[Thm]{Lemma}

\newtheorem{fig}[Thm]{Figure}

\newcommand{\kn}{\Large\bf
$K\hspace{-.4cm} N$
\Large\bf\vv }

\def\cal#1{\mathcal #1}
\def\con{{\cal O}_{n}}

\def\edot{=1,\ldots,n}
\def\pr{{\it Proof.}\quad}

\def\co#1{{\cal O}_{#1}}
\def\ltn{\ell_{2}({\bf N})}

\def\disp#1{{\displaystyle #1}}
\def\brl{branching law}
\def\bfsnl{{\rm BFS}_{N}(\Lambda)}

\addtocounter{footnote}{1}
\def\sftt#1{
\setcounter{equation}{0}
\addtocounter{footnote}{1}
\section{#1}
}

\def\ssft#1{\subsection{#1}}

\begin{document}
%
%
\def\autherp{Katsunori Kawamura}
\def\emailp{e-mail: kawamura@kurims.kyoto-u.ac.jp.}
\def\addressp{{\small {\it College of Science and Engineering, 
Ritsumeikan University,}}\\
{\small {\it 1-1-1 Noji Higashi, Kusatsu, Shiga 525-8577, Japan}}
}

\def\infw{\Lambda^{\frac{\infty}{2}}V}
\def\zhalfs{{\bf Z}+\frac{1}{2}}
\def\ems{\emptyset}
\def\pmvac{|{\rm vac}\!\!>\!\! _{\pm}}
\def\vac{|{\rm vac}\rangle _{+}}
\def\dvac{|{\rm vac}\rangle _{-}}
\def\ovac{|0\rangle}
\def\tovac{|\tilde{0}\rangle}
\def\expt#1{\langle #1\rangle}
\def\zph{{\bf Z}_{+/2}}
\def\zmh{{\bf Z}_{-/2}}
\def\brl{branching law}
\def\bfsnl{{\rm BFS}_{N}(\Lambda)}
\def\scm#1{S({\bf C}^{N})^{\otimes #1}}
\def\mqb{\{(M_{i},q_{i},B_{i})\}_{i=1}^{N}}
\def\zhalf{\mbox{${\bf Z}+\frac{1}{2}$}}
\def\zmha{\mbox{${\bf Z}_{\leq 0}-\frac{1}{2}$}}
\newcommand{\mline}{\noindent
\thicklines
\setlength{\unitlength}{.1mm}
\begin{picture}(1000,5)
\put(0,0){\line(1,0){1250}}
\end{picture}
\par
 }
\def\ptimes{\otimes_{\varphi}}
\def\qtimes{\otimes_{\tilde{\varphi}}}
\def\delp{\Delta_{\varphi}}
\def\delps{\Delta_{\varphi^{*}}}
\def\gamp{\Gamma_{\varphi}}
\def\gamps{\Gamma_{\varphi^{*}}}
\def\sem{{\sf M}}
\def\sen{{\sf N}}
\def\hdelp{\hat{\Delta}_{\varphi}}
\def\tilco#1{\tilde{\co{#1}}}
\def\ndm#1{{\bf M}_{#1}(\{0,1\})}
\def\fs{{\cal F}{\cal S}({\bf N})}
\def\ba{\mbox{\boldmath$a$}}
\def\bb{\mbox{\boldmath$b$}}
\def\bc{\mbox{\boldmath$c$}}
\def\be{\mbox{\boldmath$e$}}
\def\bp{\mbox{\boldmath$p$}}
\def\bq{\mbox{\boldmath$q$}}
\def\bu{\mbox{\boldmath$u$}}
\def\bv{\mbox{\boldmath$v$}}
\def\bw{\mbox{\boldmath$w$}}
\def\bx{\mbox{\boldmath$x$}}
\def\by{\mbox{\boldmath$y$}}
\def\bz{\mbox{\boldmath$z$}}
\def\titlepage{

\noindent
{\bf 
\noindent
\thicklines
\setlength{\unitlength}{.1mm}
\begin{picture}(1000,0)(0,-300)
\put(0,0){\kn \knn\, for \journal\, Ver.\version}
\put(0,-50){\today}
\end{picture}
}
\vspace{-2.3cm}
\quad\\
{\small file: \textsf{tit01.txt,\, J1.tex}
 \footnote{
  ${\displaystyle
  \mbox{directory: \textsf{\fileplace}, 
  file: \textsf{\incfile},\, from \startdate}}$
          }
}
\quad\\
\framebox{
 \begin{tabular}{ll}
 \textsf{Title:} &
 \begin{minipage}[t]{4in}
 \titlep
 \end{minipage}
 \\
 \textsf{Author:} &\autherp
 \end{tabular}
}
{\footnotesize	\tableofcontents }
}

%
%
%
\setcounter{section}{0}
\setcounter{footnote}{0}
\setcounter{page}{1}
\pagestyle{plain}

%
%
\title{\titlep}
\author{\autherp\thanks{\emailp}
\\
\addressp}
\date{}
\maketitle
%
%
\begin{abstract}
A new construction method of $R$-matrix is given.
Let $A$ be a C$^{*}$-bialgebra with a comultiplication $\Delta$.
For two states $\omega$ and $\psi$ of $A$
which satisfy certain conditions,
we construct a unitary $R$-matrix $R(\omega,\psi)$ of  
the C$^{*}$-bialgebra $(A,\Delta)$
on the tensor product of GNS representation spaces associated 
with $\omega$ and $\psi$.
The set  $\{R(\omega,\psi):\omega,\psi\}$ 
satisfies a kind of Yang-Baxter equation.
Furthermore, we show a nontrivial example of such $R$-matrices
for a non-quasi-cocommutative C$^{*}$-bialgebra.
\end{abstract}


%
%
\sftt{Introduction}
\label{section:first}
A C$^{*}$-bialgebra is a generalization of bialgebra in the theory of C$^{*}$-algebras,
which was introduced in C$^{*}$-algebraic framework for quantum groups \cite{KV,MNW}. 
The purpose of this paper is to show 
a new construction method of $R$-matrix
by using Gel'fand-Na\u{\i}mark-Segal (= GNS)
 representations of a C$^{*}$-bialgebra
without the assumption of the quasi-cocommutativity. 
In this section, we show our motivation,
definitions of C$^{*}$-bialgebras and the main theorem.

%
%
\ssft{Motivation}
\label{subsection:firstone}
In this subsection, we roughly explain our motivation 
and the background of this study.
Explicit mathematical definitions will 
be shown after $\S$ \ref{subsection:firsttwo}.

Let $A$ be a bialgebra with the comultiplication $\Delta$.
Recall that $(A,\Delta)$ 
is {\it quasi-cocommutative} if 
there exists an invertible element $R\in A\otimes A$
such that 
%
%
\begin{equation}
\label{eqn:rrone}
R\Delta(x) R^{-1}=\Delta^{op}(x)\quad(x\in A).
\end{equation}
Such an element $R$ is called the {\it universal $R$-matrix} of $(A,\Delta)$.
It is known that $R$ satisfies the following relation in $A\otimes A\otimes A$:
%
%
\begin{equation}
\label{eqn:ybetwo}
R_{12}R_{13}R_{23}=R_{23}R_{13}R_{12}
\end{equation}
where $R_{12}\equiv R\otimes 1$ and so on.
This is called the {\it Yang-Baxter equation}  \cite{Kassel}
(or {\it quantum Yang-Baxter equation} \cite{CP}) for $R$.
When $(A,\Delta)$ is a C$^{*}$-bialgebra,
the definition of $R$ is slightly modified 
(which will be given in Definition \ref{defi:cocommutative}(i).)
In the theory of quantum groups,
$R$-matrix and the Yang-Baxter equation are fruitful subjects
as relations with mathematical physics 
and topology \cite{Drinfeld,Jimbo,Jimbo2,Kassel,KS,KBI,Xu}
for quasi-cocommutative bialgebras.
They were also considered by Van Daele and Van Keer \cite{VanDaeleVanKeer}
for Hopf $*$-algebras.
We constructed a non-cocommutative C$^{*}$-bialgebra 
with a universal $R$-matrix,
which is defined as the direct sum of all matrix algebras in \cite{TS21}.
Furthermore, we proved that
an inductive limit of quasi-cocommutative C$^{*}$-bialgebras
is not always quasi-cocommutative in \cite{TS24}.

On the other hand,
we introduced the C$^{*}$-bialgebra $\co{*}$
as the direct sum of all Cuntz algebras except $\co{\infty}$ \cite{TS02}:
%
%
\begin{equation}
\label{eqn:cuntztwo}
\co{*}=\co{1}\oplus\co{2}\oplus\co{3}\oplus\co{4}\oplus\cdots,
\end{equation}
where $\co{1}$ denotes the $1$-dimensional C$^{*}$-algebra ${\bf C}$
for convenience.
We constructed a non-cocommutative comultiplication $\Delta_{\varphi}$ 
of ${\cal O}_{*}$. 
Unfortunately, 
there exists no universal $R$-matrix of $(\co{*},\delp)$ \cite{TS20}.

As a construction method of $R$-matrix,
the quantum double  is well-known  \cite{Drinfeld, Kassel,KS}.
In this paper, we show a new method to
 construct $R$-matrices from states of a C$^{*}$-bialgebra
under some conditions.
For two states $\omega,\psi$ of a C$^{*}$-bialgebra $A$
which satisfy some assumptions,
let ${\cal H}_{\omega}$ and ${\cal H}_{\psi}$ denote
GNS representation spaces by $\omega$ and $\psi$, respectively.
Then there exists a unitary operator $R(\omega,\psi)$
on ${\cal H}_{\omega}\otimes {\cal H}_{\psi}$ 
such that
the set $\{R(\omega,\psi):\omega,\psi\in {\cal S}\}$ 
of unitary operators satisfy 
the Yang-Baxter equation depending on a special set ${\cal S}$ of states
of $A$:
%
%
\begin{equation}
\label{eqn:ybethree}
R_{12}(\omega_{1},\omega_{2})
R_{13}(\omega_{1},\omega_{3})
R_{23}(\omega_{2},\omega_{3})
=
R_{23}(\omega_{2},\omega_{3})
R_{13}(\omega_{1},\omega_{3})
R_{12}(\omega_{1},\omega_{2}).
\end{equation}
%
Furthermore, we show a non-trivial example of this construction.

%
%
\ssft{Local $R$-matrix of C$^{*}$-bialgebra}
\label{subsection:firsttwo}
In this subsection, we recall definitions of C$^{*}$-bialgebra,
and introduce 
local $R$-matrix of a C$^{*}$-bialgebra.
At first,
we prepare terminologies about C$^{*}$-bialgebra according to \cite{KV,MNW}.
For two C$^{*}$-algebras $A$ and $B$,
let ${\rm Hom}(A,B)$ and $A\otimes B$ 
denote the set of all $*$-homomorphisms from $A$ to $B$ and 
the minimal C$^{*}$-tensor product of $A$ and $B$, respectively.
Let $M(A)$ denote the multiplier algebra of a C$^{*}$-algebra $A$.
We state that $f\in {\rm Hom}(A,M(B))$ is {\it nondegenerate} if $f (A)B$ 
is dense in a C$^{*}$-algebra $B$.
A pair $(A,\Delta)$ is a {\it C$^{*}$-bialgebra}
if $A$ is a C$^{*}$-algebra and $\Delta\in {\rm Hom}(A,M(A\otimes A))$ 
such that $\Delta$ is nondegenerate and the following holds:
%
%
\begin{equation}
\label{eqn:bialgebratwo}
(\Delta\otimes id)\circ \Delta=(id\otimes\Delta)\circ \Delta.
\end{equation}
We call $\Delta$ the {\it comultiplication} of $A$.
Remark that $A$ has no unit in general for a C$^{*}$-bialgebra $(A,\Delta)$.
Define the {\it extended flip} $\tilde{\tau}_{A,A}$ from $M(A\otimes A)$ 
to $M(A\otimes A)$ as
$\tilde{\tau}_{A,A}(X)(x\otimes y)\equiv \tau_{A,A}(X(y\otimes x))$
for $X\in M(A\otimes A),\,x, y\in A$
where $\tau_{A,A}$ denotes the flip of $A\otimes A$.
The map $\Delta^{op}$ from $A$ to $M(A\otimes A)$
defined as $\Delta^{op}\equiv \tilde{\tau}_{A,A}\circ \Delta$
is called the {\it opposite comultiplication} of $\Delta$.
A C$^{*}$-bialgebra $(A,\Delta)$ is {\it cocommutative}
if $\Delta=\Delta^{op}$. 

According to \cite{Drinfeld,Kassel,VanDaeleVanKeer},
we introduce unitary $R$-matrix, the quasi-cocommutativity
and local $R$-matrix
for C$^{*}$-bialgebra as follows.
%
%
\begin{defi}
\label{defi:cocommutative}
\begin{enumerate}
\item
An element $R$ in $M(A\otimes A)$
is called a (unitary) universal $R$-matrix of $(A,\Delta)$
if $R$ is a unitary and 
%
%
\begin{equation}
\label{eqn:univ}
R\Delta(x)R^{*}=\Delta^{op}(x)\quad
(x\in A).
\end{equation}
In this case,
$(A,\Delta)$ is said to be quasi-cocommutative
(or almost cocommutative \cite{CP}).
\item
Let $({\cal H}_{i},\pi_{i})$ be a representation of $A$
for $i=1,2$.
A unitary operator $R_{\pi_{1},\pi_{2}}$ on ${\cal H}_{1}\otimes {\cal H}_{2}$
is called a local $R$-matrix of $(A,\Delta)$ on ${\cal H}_{1}\otimes {\cal H}_{2}$
if it satisfies
%
%
\begin{equation}
\label{eqn:oppositeb}
R_{\pi_{1},\pi_{2}}(\pi_{1}\otimes \pi_{2})(\Delta(x))R_{\pi_{1},\pi_{2}}^{*}
=(\pi_{1}\otimes \pi_{2})(\Delta^{op}(x))
\quad(x\in A).
\end{equation}
%
%
%
%
\end{enumerate}
\end{defi}

\noindent
In Definition \ref{defi:cocommutative}(ii),
if $\pi_{1}=\pi_{2}$,
then the flip of ${\cal H}_{1}\otimes {\cal H}_{2}$
is a (trivial) local $R$-matrix of $(A,\Delta)$ on ${\cal H}_{1}\otimes {\cal H}_{2}$.
If $\pi_{1}\ne \pi_{2}$,
then a local $R$-matrix does not always exist when
$(A,\Delta)$ is not cocommutative.

In addition,
we introduce a new notion for C$^{*}$-bialgebra.
Let $A\odot B$ denote the algebraic tensor product 
of $*$-algebras $A$ and $B$.
%
%
\begin{defi}
\label{defi:algebraic}
A C$^{*}$-bialgebra $(A,\Delta)$ is algebraic
if there exists a dense $*$-subalgebra $A_{0}$ of the C$^{*}$-algebra $A$
such that $\Delta(A_{0})\subset A_{0}\odot A_{0}$.
We call $A_{0}$ an algebraic part of $(A,\Delta)$.
\end{defi}

\noindent
Clearly, any finite dimensional C$^{*}$-bialgebra
is algebraic.
If a C$^{*}$-bialgebra $(A,\Delta)$ is algebraic with an algebraic part $A_{0}$,
then 
$\Delta(x)$ is written as follows for any element $x\in A_{0}$:
There exist $1\leq m<\infty$ and 
$x_{1}^{'},\ldots,x_{m}^{'},x_{1}^{''},\ldots,x_{m}^{''}$ in $A_{0}$
such that 
%
%
\begin{equation}
\label{eqn:finite}
\Delta(x)=x_{1}^{'}\otimes x_{1}^{''}+\cdots +x_{m}^{'}\otimes x_{m}^{''}.
\end{equation}

%
%
\ssft{Main theorem}
\label{subsection:firstthree}
In this subsection, we show our main theorem.
For a C$^{*}$-algebra $A$,
let ${\cal S}(A)$ and ${\rm Rep}A$
denote the set of all states 
and the class of all nondegenerate representations of $A$, respectively.
%
%
\begin{defi}
\label{defi:gnsmap}
\cite{ES,KV,MNW}
For $\omega\in {\cal S}(A)$ of a C$^{*}$-algebra $A$
with the GNS triple 
$({\cal H}_{\omega},\pi_{\omega},\Omega_{\omega})$,
the linear map $\Lambda_{\omega}$ from $A$ to ${\cal H}_{\omega}$
defined by
%
%
\begin{equation}
\label{eqn:gnstwo}
\Lambda_{\omega}(x)\equiv \pi_{\omega}(x)\Omega_{\omega}\quad(x\in A),
\end{equation}
is called the GNS map associated with $\omega$. 
\end{defi}
Then our main theorem is stated as follows.
%
%
\begin{Thm}
\label{Thm:maintwo}
Let $(A,\Delta)$ be an algebraic C$^{*}$-bialgebra.
For $\pi_{1},\pi_{2}\in {\rm Rep}A$ 
and $\omega_{1},\omega_{2}\in {\cal S}(A)$,
we write
%
%
\begin{equation}
\label{eqn:starone}
\pi_{1}\star \pi_{2}\equiv (\pi_{1}\otimes \pi_{2})\circ \Delta,\quad
\omega_{1}\star \omega_{2}\equiv (\omega_{1}\otimes \omega_{2})\circ \Delta
\end{equation}
where $(\omega_{1}\otimes \omega_{2})(x\otimes y)
\equiv \omega_{1}(x)\omega_{2}(y)$ for $x,y\in A$ (\cite{KR2}, p.847).
Assume that a non-empty subset ${\cal S}_{0}$ of ${\cal S}(A)$ 
satisfies the following conditions:
\def\labelenumi{\theenumi}
\def\theenumi{{\rm (\alph{enumi})}}
\begin{enumerate}
\item
The set ${\cal S}_{0}$ is closed with respect to the operation $\star$
in (\ref{eqn:starone}) and the semigroup $({\cal S}_{0},\star)$ is abelian.
\item
For any $\omega_{1},\omega_{2}\in {\cal S}_{0}$,
$\Omega_{\omega_{1}}\otimes \Omega_{\omega_{2}}$ is a cyclic 
vector for 
the representation $({\cal H}_{\omega_{1}}\otimes {\cal H}_{\omega_{2}},\,
\pi_{\omega_{1}}\star \pi_{\omega_{2}})$ of $A$.
\end{enumerate}
\def\labelenumi{\theenumi}
\def\theenumi{{\rm (\roman{enumi})}}
Define the unitary $R(\omega_{1},\omega_{2})$ from
${\cal H}_{\omega_{1}}\otimes {\cal H}_{\omega_{2}}$
to ${\cal H}_{\omega_{1}}\otimes {\cal H}_{\omega_{2}}$ by
%
%
\begin{equation}
\label{eqn:rmattwo}
R(\omega_{1},\omega_{2})\Lambda_{\omega_{1},\omega_{2}}(\Delta(x))
\equiv \Lambda_{\omega_{1},\omega_{2}}(\Delta^{op}(x))\quad(x\in A)
\end{equation}
where $\Lambda_{\omega_{1},\omega_{2}}
\equiv \Lambda_{\omega_{1}}\otimes\Lambda_{\omega_{2}}$.
Then $R(\omega_{1},\omega_{2})$ is well-defined and the following holds:
\begin{enumerate}
\item
For any $\omega_{1},\omega_{2}\in {\cal S}_{0}$
and $x\in A$,
the following identities hold:
%
%
\begin{eqnarray}
\label{eqn:implement}
\hspace{-1cm}R(\omega_{1},\omega_{2})
(\pi_{\omega_{1}}\otimes \pi_{\omega_{2}})(\Delta(x))
(R(\omega_{1},\omega_{2}))^{*}
=&(\pi_{\omega_{1}}\otimes \pi_{\omega_{2}})(\Delta^{op}(x)),\\
\nonumber
\\
\label{eqn:symone}
R(\omega_{1},\omega_{2})\tau_{\omega_{2},\omega_{1}}R(\omega_{2},\omega_{1})
\tau_{\omega_{1},\omega_{2}}=& I_{\omega_{1}}\otimes I_{\omega_{2}}
\end{eqnarray}
where $\tau_{\omega_{1},\omega_{2}}$ denotes
the flip of ${\cal H}_{\omega_{1}}\otimes {\cal H}_{\omega_{2}}$.
\item
For any $\omega_{1},\omega_{2}\in {\cal S}_{0}$,
$\pi_{\omega_{1}}\star \pi_{\omega_{2}}$
and 
$\pi_{\omega_{2}}\star \pi_{\omega_{1}}$
are unitarily equivalent.
\item
For any $\omega_{1},\omega_{2},\omega_{3}\in {\cal S}_{0}$,
the following identity holds:
%
%
\begin{equation}
\label{eqn:ybeone}
R_{12}(\omega_{1},\omega_{2})
R_{13}(\omega_{1},\omega_{3})
R_{23}(\omega_{2},\omega_{3})
=
R_{23}(\omega_{2},\omega_{3})
R_{13}(\omega_{1},\omega_{3})
R_{12}(\omega_{1},\omega_{2})
\end{equation}
on  ${\cal H}_{\omega_{1}}\otimes {\cal H}_{\omega_{2}}
\otimes {\cal H}_{\omega_{3}}$
where we use the leg numbering notation in \cite{BS}.
For example,
$R_{12}(\omega_{1},\omega_{2})
\equiv R(\omega_{1},\omega_{2})\otimes I_{\omega_{3}} $.
\end{enumerate}
\end{Thm}

\noindent
From (\ref{eqn:implement}),
$R(\omega_{1},\omega_{2})$ is a local $R$-matrix of 
$(A,\Delta)$ on ${\cal H}_{\omega_{1}}\otimes {\cal H}_{\omega_{2}}$.
The equation (\ref{eqn:ybeone}) is regarded as a kind of Yang-Baxter equation.

%
%
\begin{rem}
\label{rem:one}
{\rm
\begin{enumerate}
\item
If $\omega_{1}=\omega_{2}$,
then $R(\omega_{1},\omega_{2})=I_{\omega_{1}}\otimes I_{\omega_{2}}$.
If not,
then  this does not hold in general.
If $({\cal S}_{0},\star)$ is not abelian,
then (\ref{eqn:implement}) does not hold in general.
We will  show these examples in $\S$ \ref{subsection:thirdthree}.
\item
We explain conditions in Theorem \ref{Thm:maintwo}. 
The condition (b)
does not always hold even if the condition (a) holds.
For $\pi\in {\rm Rep}A$,
let $[\pi]$ denote the unitary equivalence class of $\pi$.
For $\pi_{1},\pi_{2}\in {\rm Rep}A$,
the new product
$[\pi_{1}]\star [\pi_{2}]\equiv [\pi_{1}\star \pi_{2}]$ is well-defined.
Let ${\sf R}(A)$ denote the set of unitary equivalence classes
of all representations of $A$.
Then we obtain the map
$f$ from ${\cal S}_{0}$ to ${\sf R}(A)$ by
$f(\omega)\equiv [\pi_{\omega}]$ for $\omega\in {\cal S}_{0}$.
Then the condition in Theorem \ref{Thm:maintwo}
is interpreted as follows:
%
%
\begin{equation}
\label{eqn:f}
f(\omega\star\psi)=f(\omega)\star f(\psi)\quad(\omega,\psi\in {\cal S}_{0}).
\end{equation}
That is,
$f$ is a semigroup homomorphism from 
$({\cal S}_{0},\star)$ to $({\sf R}(A),\star)$.
Such examples will be shown in $\S$ \ref{subsection:thirdtwo}.
\item
The GNS map in Definition \ref{defi:gnsmap} was used in 
constructions of Kac-Takesaki operators \cite{ES,KV,MNW}
and an analogue of multiplicative isometry \cite{TS08}.
\end{enumerate}
}
\end{rem}

In $\S$ \ref{section:second},
we will prove Theorem \ref{Thm:maintwo}.
In $\S$ \ref{section:third},
we will treat $(\co{*},\delp)$ in (\ref{eqn:cuntztwo})
and certain states of $\co{*}$
as an example of Theorem \ref{Thm:maintwo}.

%
%
\sftt{Proof of main theorem}
\label{section:second}
In this section,
we prove Theorem \ref{Thm:maintwo}.
%
%
\ssft{Lemma for GNS maps}
\label{subsetion:secondone}
\label{subsection:secondone}
In order to prove Theorem \ref{Thm:maintwo},
we prepare an elementary lemma for GNS maps in (\ref{eqn:gnstwo})
in this subsection.
%
%
\begin{lem}
\label{lem:gns}
Let $\Lambda_{\omega}$ be as in (\ref{eqn:gnstwo})
and let $A,B,C$ be C$^{*}$-algebras.
\begin{enumerate}
\item
For $\phi\in {\rm Hom}(A,B)$
and a state $\omega$ of $B$,
define
%
%
\begin{equation}
\label{eqn:unitaryone}
U\Lambda_{\omega\circ \phi}(x)\equiv \Lambda_{\omega}(\phi(x))\quad(x\in A).
\end{equation}
Then $U$ is an isometry 
from ${\cal H}_{\omega\circ \phi}$ to ${\cal H}_{\omega}$
such that 
$U^{*}\pi_{\omega}(\phi(x))U=\pi_{\omega\circ \phi}(x)$ for $x\in A$.
\item
In addition to (i),
if $\Omega_{\omega}$ is a cyclic vector 
for $({\cal H}_{\omega},\pi_{\omega}\circ \phi)$,
then $U$ is a unitary.
\item
Let $\phi\in {\rm Hom}(A,B\otimes C)$ and 
let $\omega_{1}$ and $\omega_{2}$ be states of $B$ and $C$, respectively.
Let $\omega\equiv (\omega_{1}\otimes \omega_{2})\circ \phi$.
If 
$\Omega_{\omega_{1}}\otimes \Omega_{\omega_{2}}$
is a cyclic vector for 
$({\cal H}_{\omega_{1}}\otimes {\cal H}_{\omega_{2}},
(\pi_{\omega_{1}}\otimes \pi_{\omega_{2}})\circ \phi)$,
then 
%
%
\begin{equation}
\label{eqn:unitarytwo}
U(\omega_{1},\omega_{2};\phi)
\Lambda_{\omega}(x)\equiv \Lambda_{\omega_{1},\omega_{2}}(\phi(x))
\quad(x\in A)
\end{equation}
defines a unitary $U(\omega_{1},\omega_{2};\phi)$
from ${\cal H}_{\omega}$ to 
${\cal H}_{\omega_{1}}\otimes {\cal H}_{\omega_{2}}$
such that
%
%
\begin{equation}
\label{eqn:unitarythree}
(U(\omega_{1},\omega_{2};\phi))^{*}
(\pi_{\omega_{1}}\otimes \pi_{\omega_{2}})(\phi(x))
U(\omega_{1},\omega_{2};\phi)=\pi_{\omega}(x)\quad
(x\in A).
\end{equation}
%
\item
In addition to the assumption in (iii),
assume that $B=C$ and the following identity holds:
%
%
\begin{equation}
\label{eqn:commute}
(\omega_{1}\otimes \omega_{2})\circ \phi=
(\omega_{2}\otimes \omega_{1})\circ \phi.
\end{equation}
Define the unitary $R(\omega_{1},\omega_{2};\phi)$ from 
${\cal H}_{\omega_{1}}\otimes {\cal H}_{\omega_{2}}$ to 
${\cal H}_{\omega_{1}}\otimes {\cal H}_{\omega_{2}}$
by
%
%
\begin{equation}
\label{eqn:unitaryfour}
R(\omega_{1},\omega_{2};\phi)
\Lambda_{\omega_{1},\omega_{2}}(\phi(x))
\equiv 
\Lambda_{\omega_{1},\omega_{2}}(\phi^{op}(x))\quad(x\in A)
\end{equation}
where $\phi^{op}\equiv \tau\circ \phi$
and $\tau$ denotes the flip of $B\otimes B$.
Then $R(\omega_{1},\omega_{2};\phi)$ is well-defined and
the following holds:
%
%
\begin{equation}
\label{eqn:unitaryfive}
R(\omega_{1},\omega_{2};\phi)(\pi_{\omega_{1}}\otimes 
\pi_{\omega_{2}})(\phi(x))(R(\omega_{1},\omega_{2};\phi))^{*}
=(\pi_{\omega_{1}}\otimes \pi_{\omega_{2}})(\phi^{op}(x))
\end{equation}
for $x\in A$.
\end{enumerate}
\end{lem}
%
%
\pr
From the uniqueness of the GNS representation,
(i) and (ii) hold. 
For example,  see Proposition 4.5.3 in \cite{KR1}.

\noindent
(iii)
The statement holds from (i) and (ii).

\noindent
(iv)
By (\ref{eqn:commute}), we obtain
${\cal H}_{(\omega_{1}\otimes \omega_{2})\circ \phi^{op}}
={\cal H}_{(\omega_{2}\otimes \omega_{1})\circ \phi}
={\cal H}_{(\omega_{1}\otimes \omega_{2})\circ \phi}$.
From this and the definition of $R(\omega_{1},\omega_{2};\phi)$,
we see that
%
%
\begin{equation}
\label{eqn:equal}
R(\omega_{1},\omega_{2};\phi)=
U(\omega_{1}, \omega_{2};\phi^{op})
(U(\omega_{1},\omega_{2};\phi))^{*}.
\end{equation}
Hence the statement holds from (iii).
\qedh

\noindent
We illustrate Lemma \ref{lem:gns}(iv) as follows:

\def\diagram{
\put(0,0){${\cal H}_{\omega_{1}}\otimes {\cal H}_{\omega_{2}}$}
\put(800,0){${\cal H}_{\omega_{1}}\otimes {\cal H}_{\omega_{2}}$}
\put(200,-200){
${\cal H}_{(\omega_{1}\otimes \omega_{2})\circ \phi}\quad =\quad
{\cal H}_{(\omega_{1}\otimes \omega_{2})\circ \phi^{op}}$}
\put(250,10){\vector(1,0){470}}
\put(140,-50){\vector(1,-1){100}}
\put(700,-150){\vector(1,1){100}}
\put(360,40){$R(\omega_{1},\omega_{2}:\phi)$}
\put(-110,-130){$(U(\omega_{1},\omega_{2};\phi))^{*}$}
\put(800,-130){$U(\omega_{1},\omega_{2};\phi^{op})$}
\put(430,-100){\scalebox{2}[2]{$\circlearrowleft$}}
}
\noindent
\thicklines
%
%
\begin{fig}
\label{fig:one}
\quad\\
\setlength{\unitlength}{.1mm}
\begin{picture}(1200,300)(-140,-230)
\put(0,0){\diagram}
\end{picture}
\end{fig}

%
%
\begin{rem}
\label{rem:commutative}
{\rm 
Lemma \ref{lem:gns} holds
without the assumption that $\phi$ is a comultiplication.
Since the logic of the proof is very elementary, 
the assumption about states is essential.
}
\end{rem}

%
%
\ssft{Proof of Theorem \ref{Thm:maintwo}}
\label{subsection:secondtwo}
We prove Theorem \ref{Thm:maintwo} in this subsection.
Let $(A,\Delta)$ be as in Theorem \ref{Thm:maintwo}.
Applying Lemma \ref{lem:gns}(iv)
to the case $B=C=A$ and $\phi=\Delta$,
we see that $R(\omega_{1},\omega_{2})=R(\omega_{1},\omega_{2};\Delta)$.
Hence $R(\omega_{1},\omega_{2})$ 
in (\ref{eqn:rmattwo})
is well-defined.

\noindent
(i)
From Lemma \ref{lem:gns}(iv), (\ref{eqn:implement}) holds.
We can verify that 
(\ref{eqn:symone}) holds on $\Lambda_{\omega_{1},\omega_{2}}(\Delta^{op}(x))$
for each $x\in A$.
Hence (\ref{eqn:symone}) holds
because
$\{\Lambda_{\omega_{1},\omega_{2}}(\Delta^{op}(x)):x\in A\}$ 
is dense in ${\cal H}_{\omega_{1}}\otimes {\cal H}_{\omega_{2}}$.

\noindent
(ii)
This follows from (i).

\noindent
(iii)
Define
%
%
\begin{equation}
\label{eqn:proof}
F_{R}\equiv (id\otimes \Delta)\circ \Delta,\quad
F_{L}\equiv (\Delta\otimes id)\circ \Delta.
\end{equation}
From (\ref{eqn:bialgebratwo}), $F_{R}=F_{L}$.
Let $A_{0}$ be an algebraic part of $(A,\Delta)$.
By assumption,
we can write as follows for $x\in A_{0}$:
%
%
\begin{eqnarray}
\label{eqn:pt}
F_{R}(x)=&
\disp{\sum_{i} x_{i}^{'}\otimes (x_{i}^{''})^{'}\otimes (x_{i}^{''})^{''},}\\
\nonumber
\\
\label{eqn:ps}
F_{L}(x)=&
\disp{\sum_{j} (y_{j}^{'})^{'}\otimes (y_{j}^{'})^{''}\otimes y_{j}^{''}}
\end{eqnarray}
where R.H.S.s of  
(\ref{eqn:pt}) and (\ref{eqn:ps}) are finite sums
of elements in $A_{0}$.

Here we write elements in ${\cal S}_{0}$
as $a,b,c,\ldots$ for simplicity of description.
From (\ref{eqn:unitarythree}) and the assumption
for $\Omega_{a}\otimes\Omega_{b}$, 
$\pi_{a}\star \pi_{b}$ is unitarily equivalent to $\pi_{a\star b}$
for each $a,b\in {\cal S}_{0}$.
From this,
$(\pi_{a}\star \pi_{b})\star\pi_{c}$ is unitarily equivalent to $\pi_{a\star b}\star\pi_{c}$.
By assumption,  
$\Omega_{a\star b}\otimes \Omega_{c}$ is a cyclic vector for 
$({\cal H}_{a\star b}\otimes {\cal H}_{b},\pi_{a\star b}\star\pi_{c})$. 
Hence 
$\Omega_{a}\otimes \Omega_{b}\otimes \Omega_{c}$ is a also cyclic vector for 
$(\pi_{a}\star\pi_{b})\star\pi_{c}$.
From this,
$\Omega_{a}\otimes \Omega_{b}\otimes \Omega_{c}$ is a cyclic vector for 
$(\pi_{a}\otimes \pi_{b}\otimes \pi_{c})\circ F_{R}$
($=
(\pi_{a}\otimes \pi_{b}\otimes \pi_{c})\circ F_{L}$)  on 
${\cal H}_{a}\otimes {\cal H}_{b}\otimes {\cal H}_{c}$.
Therefore
$\{\Lambda_{a,b,c}(F_{R}(x)):x\in A_{0}\}$
is dense in ${\cal H}_{a}\otimes {\cal H}_{b}\otimes {\cal H}_{c}$
where $\Lambda_{a,b,c}\equiv \Lambda_{a}\otimes\Lambda_{b}\otimes\Lambda_{c}$.
Hence it is sufficient to show 
(\ref{eqn:ybeone}) on $\{\Lambda_{a,b,c}(F_{R}(x)):x\in A_{0}\}$.
The following holds:\\
\\
$R_{12}(a,b)
R_{13}(a,c)
R_{23}(b,c)
\Lambda_{a,b,c}(F_{R}(x))$
\[
\begin{array}{rl}
=&
R_{13}(a,c)
R_{23}(b,c)
\Lambda_{a,b,c}((id\otimes \Delta^{op})(\Delta(x)))\\
\\
=&
\disp{
\sum R_{12}(a,b)
R_{13}(a,c)
\Lambda_{a,b,c}
(\,x_{i}^{'}\otimes (x_{i}^{''})^{''}\otimes (x_{i}^{'})^{'}\,)}\\
\\
=&
\disp{
\sum R_{12}(a,b)
R_{13}(a,c)
\Lambda_{a,b,c}
(\,(y_{j}^{'})^{'}\otimes y_{j}^{''}\otimes (y_{j}^{'})^{''}\,)}
\quad(\mbox{by }F_{R}(x)=F_{L}(x))\\
\\
=&
\disp{
\sum R_{12}(a,b)
\Lambda_{a,b,c}(\,(y_{j}^{'})^{''}\otimes y_{j}^{''}\otimes (y_{j}^{'})^{'}\,)}\\
\\
=&
\disp{
\sum R_{12}(a,b)
\Lambda_{a,b,c}
(\,(x_{i}^{'})^{'}\otimes (x_{i}^{''})^{''}\otimes x_{i}^{'}\,)}
\quad(\mbox{by }F_{R}(x)=F_{L}(x))\\
\\
=&
\disp{
\sum 
\Lambda_{a,b,c}
(\,(x_{i}^{'})^{''}\otimes (x_{i}^{'})^{''}\otimes x_{i}^{'}\,)}\\
\\
=&
\disp{
\Lambda_{a,b,c}((\Delta^{op}\otimes id)(\Delta^{op}(x))).}\\
\end{array}
\]
Remark that every ``$\Sigma$" means  a finite sum.
As the same token,
%
%
\begin{equation}
\label{eqn:lefttwo}
R_{23}(b,c)
R_{13}(a,c)
R_{12}(a,b)
\Lambda_{a,b,c}(F_{L}(x))
=
\Lambda_{a,b,c}((id\otimes \Delta^{op})(\Delta^{op}(x))).
\end{equation}
Hence the statement holds
from the  coassociativity of $\Delta^{op}$.
\qedh

%
%
\sftt{C$^{*}$-bialgebra defined as the direct sum of Cuntz algebras
and its local $R$-matrices}
\label{section:third}
We give a set of states 
of the C$^{*}$-algebra $\co{*}$ in (\ref{eqn:cuntztwo})
which satisfies assumptions in Theorem \ref{Thm:maintwo} in this section.
In $\S$ \ref{subsection:thirdthree},
we will show a concrete local $R$-matrix.
%
%
\ssft{C$^{*}$-bialgebra $(\co{*},\delp)$}
\label{subsection:thirdone}
In this subsection, we recall the C$^{*}$-bialgebra $(\co{*},\delp)$ in \cite{TS02}.
For $n\geq 2$, 
consider  the {\it Cuntz algebra} $\con$  \cite{C}, that is,
a C$^{*}$-algebra which is universally generated by
generators $s_{1},\ldots,s_{n}$ satisfying
$s_{i}^{*}s_{j}=\delta_{ij}I$ for $i,j\edot$ and
$\sum_{i=1}^{n}s_{i}s_{i}^{*}=I$
where $I$ denotes the unit of $\con$.
Since $\con$ is simple, that is, there is no
non-trivial closed two-sided ideal,
any $*$-homomorphism from $\con$ to a C$^{*}$-algebra is injective.
If $t_{1},\ldots,t_{n}$ are elements of a unital C$^{*}$-algebra $A$ such that
$t_{1},\ldots,t_{n}$ satisfy the relations of canonical generators of $\con$,
then the correspondence $s_{i}\mapsto t_{i}$ for $i\edot$
is uniquely extended to a $*$-embedding
of $\con$ into $A$ from the uniqueness of $\con$.
Therefore we simply call such a correspondence 
among generators by an embedding of $\con$ into $A$.

Redefine the C$^{*}$-algebra $\co{*}$ in (\ref{eqn:cuntztwo})
as the direct sum of the set $\{\con:n\in {\bf N}\}$ of Cuntz algebras:
%
%
\begin{equation}
\label{eqn:cuntbi}
\co{*}\equiv \bigoplus_{n\in {\bf N}} \con
=\{(x_{n}):\|x_{n}\|\to 0\mbox{ as }n\to\infty\}
\end{equation}
where ${\bf N}=\{1,2,3,\ldots\}$
and $\co{1}$ denotes the $1$-dimensional C$^{*}$-algebra for convenience.
For $n\in {\bf N}$, let $I_{n}$ denote the unit of $\con$ 
and let $s_{1}^{(n)},\ldots,s_{n}^{(n)}$ denote
canonical generators of $\con$
where $s_{1}^{(1)}\equiv I_{1}$.
For $n,m\in {\bf N}$,
define $\varphi_{n,m}\in {\rm Hom}(\co{nm},\con\otimes \co{m})$ by
%
%
\begin{equation}
\label{eqn:embeddingone}
\varphi_{n,m}(s_{m(i-1)+j}^{(nm)})\equiv s_{i}^{(n)}\otimes s_{j}^{(m)}
\quad(i=1,\ldots,n,\,j=1,\ldots,m).
\end{equation}
%
%
\begin{Thm}
(\cite{TS02}, Theorem 1.1)
\label{Thm:mainone}
For the set $\varphi\equiv \{\varphi_{n,m}:n,m\in {\bf N}\}$ in
(\ref{eqn:embeddingone}),
define the $*$-homomorphism $\delp$ from $\co{*}$ to $\co{*}\otimes \co{*}$ by
%
%
\begin{eqnarray}
\label{eqn:dpone}
\delp\equiv& \oplus\{\delp^{(n)}:n\in {\bf N}\},\\
\nonumber
\\
\label{eqn:dptwo}
\delp^{(n)}(x)\equiv &\disp{\sum_{(m,l)\in {\bf N}^{2},\,ml=n}\varphi_{m,l}(x)}
\quad(x\in \con,\,n\in {\bf N}).
\end{eqnarray}
Then  $(\co{*},\delp)$ is an algebraic  C$^{*}$-bialgebra such that
$\delp(\co{*})\subset \co{*}\otimes \co{*}$.
\end{Thm}

\noindent
About properties of $(\co{*},\delp)$, see \cite{TS02,TS08}.
About a generalization of $(\co{*},\delp)$, see \cite{TS05}.

Let ${\rm Rep}\con$ denote the class of all $*$-representations of $\con$.
For $\pi_{1},\pi_{2}\in{\rm Rep}\con$,
we define the relation 
$\pi_{1}\sim \pi_{2}$ if $\pi_{1}$ and $\pi_{2}$ are unitarily equivalent.
Then the following holds.
%
%
\begin{lem}(\cite{TS01}, Lemma 1.2)
\label{lem:fundamental}
For $\varphi_{n,m}$ in (\ref{eqn:embeddingone}),
$\pi_{1}\in {\rm Rep}\con$ and  $\pi_{2}\in {\rm Rep}\co{m}$,
define $\pi_{1}\ptimes \pi_{2}\in {\rm Rep}\co{nm}$ by
%
%
\begin{equation}
\label{eqn:ptimes}
\pi_{1}\ptimes \pi_{2}\equiv (\pi_{1}\otimes \pi_{2})\circ \varphi_{n,m}.
\end{equation}
Then the following holds for
$\pi_{1},\pi_{1}^{'}\in {\rm Rep}\con$,
$\pi_{2},\pi_{2}^{'}\in {\rm Rep}\co{m}$ and $\pi_{3}\in {\rm Rep}\co{l}$:
\begin{enumerate}
\item
If $\pi_{1}\sim \pi_{1}^{'}$ and $\pi_{2}\sim \pi_{2}^{'}$,
then $\pi_{1}\otimes_{\varphi} \pi_{2}\sim
\pi_{1}^{'}\otimes_{\varphi} \pi_{2}^{'}$.
\item
$\pi_{1}\otimes_{\varphi} (\pi_{2}\oplus \pi_{2}^{'})=
\pi_{1}\otimes_{\varphi} \pi_{2}\,\oplus\, \pi_{1}\otimes_{\varphi} \pi_{2}^{'}$.
\item
$\pi_{1}\otimes_{\varphi} (\pi_{2}\otimes_{\varphi} \pi_{3})
=(\pi_{1}\otimes_{\varphi} \pi_{2})\otimes_{\varphi} \pi_{3}$.
\end{enumerate}
\end{lem}

\noindent
From Lemma \ref{lem:fundamental}(i),
we can define $[\pi_{1}]\ptimes [\pi_{2}]\equiv [\pi_{1}\ptimes \pi_{2}]$
where  $[\pi]$  denotes the unitary equivalence class of $\pi$.

Let ${\cal S}_{n}$ denote the set of all states of $\con$.
For $(\omega,\omega^{'})\in {\cal S}_{n}\times {\cal S}_{m}$,
define
%
%
\begin{equation}
\label{eqn:eleven}
\omega\ptimes \omega^{'}\equiv (\omega\otimes \omega^{'})\circ
\varphi_{n,m}
\end{equation}
where $(\omega\otimes \omega^{'})(x\otimes y)\equiv \omega(x)\omega^{'}(y)$
for $x\in \con$ and $y\in\co{m}$.
Then we see that 
$\omega \ptimes (\omega^{'}\ptimes \omega^{''})
=
(\omega \ptimes \omega^{'})\ptimes \omega^{''}$.

By identifying $\con$ with 
a C$^{*}$-subalgebra of $\co{*}$, 
any state and representation of $\con$ are naturally identified with those of $\co{*}$.
From this and the definition of $\delp$,
the following holds:
%
%
\begin{equation}
\label{eqn:identify}
\pi_{1}\ptimes \pi_{2}=(\pi_{1}\otimes \pi_{2})\circ \delp,\quad
\omega_{1}\ptimes \omega_{2}=(\omega_{1}\otimes \omega_{2})\circ \delp.
\end{equation}
%

%
%
\ssft{Pure states of Cuntz algebras parametrized by unit vectors}
\label{subsection:thirdtwo}
In this subsection, we show examples of set of 
states which satisfies assumptions in Theorem \ref{Thm:maintwo}.
We recall states in \cite{GP0123}
and show tensor product formulas among them.
Let $S({\bf C}^{n})\equiv 
\{(z_1,\ldots,z_{n})\in {\bf C}^{n}:|z_{1}|^{2}+\cdots +|z_{n}|^{2}=1\}$ 
denote the set of all unit vectors in ${\bf C}^{n}$
for $n\geq 2$.
%
%
\begin{defi}
\label{defi:firstb}(\cite{GP0123}, Proposition 3.1)
For $n\geq 2$,
let $s_{1},\ldots,s_{n}$ denote canonical generators of $\con$.
For $z=(z_{1},\ldots,z_{n})\in S({\bf C}^{n})$,
define the state $\varrho_{z}$ of $\con$ by
%
%
\begin{equation}
\label{eqn:gpstate}
\varrho_{z}(s_{j_{1}}\cdots s_{j_{a}}s_{k_{b}}^{*}\cdots s_{k_{1}}^{*})\equiv 
\overline{z}_{j_{1}}\cdots \overline{z}_{j_{a}}z_{k_{b}}\cdots z_{k_{1}}
\end{equation}
for each $j_{1},\ldots,j_{a},k_{1},\ldots,k_{b}\in \{1,\ldots,n\}$
and $a,b\geq 1$.
\end{defi}
The following results for $\varrho_{z}$ are known:
For any $z$, $\varrho_{z}$ is pure when $n\geq 2$. 
Define the state $\varrho_{1}$ of $\co{1}={\bf C}$
by $\varrho_{1}(x)=x$ for each $x\in \co{1}$.
We define $S({\bf C}^{1})\equiv \{1\}\subset {\bf C}$
for convenience.
Then $\varrho_{z}$ in (\ref{eqn:gpstate}) makes sense for each
$z\in \bigcup_{n\geq 1}S({\bf C}^{n})$.
If $z,y\in S({\bf C}^{n})$ and $z\ne y$,
then GNS representations associated with $\varrho_{z}$ and $\varrho_{y}$
are not unitarily equivalent.

Let $GP(z)$ denote
the unitary equivalence class of the GNS representation associated 
with $\varrho_{z}$.
If $\pi_{1}$ and $\pi_{2}$
are representatives of $GP(z)$ and $GP(y)$ for 
$z\in S({\bf C}^{n})$ and $y\in S({\bf C}^{m})$, respectively,
then we write $GP(z)\ptimes GP(y)$ as $[\pi_{1}]\ptimes [\pi_{2}]$
for simplicity of description. 
%
%
\begin{Thm}(\cite{TS08}, Theorem 3.2)
\label{Thm:commuteb}
For $\ptimes$ in (\ref{eqn:ptimes}),
the following holds
for each $z=(z_{i})_{i=1}^{n}\in S({\bf C}^{n})$ 
and $y=(y_{j})_{j=1}^{m}\in S({\bf C}^{m})$:
\begin{enumerate}
\item
$\varrho_{z}\ptimes \varrho_{y}=\varrho_{z\boxtimes  y}$,
\item
$GP(z)\ptimes GP(y)=GP(z\boxtimes y)$
\end{enumerate}
where $z\boxtimes y\in S({\bf C}^{nm})$ is defined as 
%
%
\begin{equation}
\label{eqn:zy}
(z\boxtimes y)_{m(i-1)+j}\equiv z_{i}y_{j}\quad(i=1,\ldots,n,\,j=1,\ldots,m).
\end{equation}
\end{Thm}

For $z\in S({\bf C}^{n})$,
let $({\cal H}_{z},\pi_{z},\Omega_{z})$ denote the GNS triple 
associated with $\varrho_{z}$.
From Theorem \ref{Thm:commuteb}(ii),
$({\cal H}_{z}\otimes {\cal H}_{y},\pi_{z}\ptimes \pi_{y})$
is irreducible for each $z,y\in \bigcup_{n\geq 1}S({\bf C}^{n})$.
From this,
$\Omega_{z}\otimes \Omega_{y}$ is a cyclic vector 
for $({\cal H}_{z}\otimes {\cal H}_{y},\pi_{z}\ptimes \pi_{y})$.
In consequence,
%
%
\begin{equation}
\label{eqn:last}
{\cal S}_{0}\equiv  \bigcup_{n\geq 1}\{\varrho_{z}:z\in S({\bf C}^{n})\}
\end{equation}
satisfies all assumptions in Theorem \ref{Thm:maintwo}
with respect to $(\co{*},\delp)$.
Furthermore, if ${\cal M}$ is a subsemigroup of 
the semigroup $(\bigcup_{n\geq 1}S({\bf C}^{n}),\,\boxtimes)$,
then
$\{\varrho_{z}:z\in {\cal M}\}$ also satisfies 
all assumptions in Theorem \ref{Thm:maintwo}.
Such subsemigroups are shown in the last paragraph of $\S$ 3 in \cite{TS08}.

%
%
\ssft{Examples of local $R$-matrices}
\label{subsection:thirdthree}
In this subsection, we show two examples of Theorem \ref{Thm:maintwo}
by using special states in $\S$ \ref{subsection:thirdtwo}.
%
%
\begin{ex}
\label{ex:third}
{\rm
We show a nontrivial example of Theorem \ref{Thm:maintwo}.
Let $\{s^{(n)}_{i}\}_{i=1}^{n}$ denote the canonical generators of $\con$
and let $\ltn$ denote the Hilbert space with the orthonormal basis
$\{e_{k}:k\in {\bf N}\}$ where ${\bf N}=\{1,2,3,\ldots\}$.
Define the representation $\pi_{n}$ of $\con$ on $\ltn$ by
%
%
\begin{equation}
\label{eqn:pin}
\pi_{n}(s_{i}^{(n)})e_{k}\equiv e_{n(k-1)+i}\quad(i=1,\ldots,n,\,k\in {\bf N}).
\end{equation}
Define the state $\omega_{n}$ of $\con$
associated with $\pi_{n}$ and the unit vector $e_{1}\in\ltn$:
%
%
\begin{equation}
\label{eqn:vectorone}
\omega_{n}=\langle e_{1}|\pi_{n}(\cdot)e_{1}\rangle \quad(n\geq 1).
\end{equation}
Then $\omega_{n}$ is $\varrho_{z}$ in (\ref{eqn:gpstate})
for $z=(1,0,\ldots,0)\in S({\bf C}^{n})$.
Then $\omega_{n}\star\omega_{m}=\omega_{nm}$ for each $n,m\geq 2$.
We identify $\pi_{n}$ with the GNS representation of $\con$ by $\omega_{n}$,
that is,
$({\cal H}_{\omega_{n}},\pi_{\omega_{n}},\Omega_{\omega_{n}})
=(\ltn,\pi_{n},e_{1})$.
For the GNS map $\Lambda_{\omega}$ in (\ref{eqn:gnstwo}),
rewrite $\Lambda_{n}\equiv \Lambda_{\omega_{n}}$
and $\Lambda_{n,m}\equiv \Lambda_{\omega_{n},\omega_{m}}$.
By definition,
$\Lambda_{n}(s_{i}^{(n)})=e_{i}$ for $i=1,\ldots,n$.
From this and (\ref{eqn:dpone}),
%
%
\begin{equation}
\label{eqn:lambdathree}
\Lambda_{n,m}(\delp(s_{m(i-1)+j}^{(nm)}))
=e_{i}\otimes e_{j}
\quad(i=1,\ldots,n,\,j=1,\ldots,m).
\end{equation}
For the local $R$-matrix $R(\omega_{1},\omega_{2})$ in (\ref{eqn:rmattwo}),
rewrite $R^{(n,m)}\equiv R(\omega_{n},\omega_{m})$.
By definition and (\ref{eqn:lambdathree}),
we can verify that
%
%
\begin{equation}
\label{eqn:newramatrix}
R^{(n,m)}(e_{i}\otimes e_{j})= e_{\underline{i}}\otimes e_{\underline{j}}
\quad ((i,j)\in
\{1,\ldots,n\}\times \{1,\ldots,m\})
\end{equation}
where
$(\underline{i},\underline{j})\in
\{1,\ldots,n\}\times \{1,\ldots,m\}$ is uniquely determined as
the following linear Diophantine equation:
%
%
\begin{equation}
\label{eqn:dio}
m(i-1)+j=n(\underline{j}-1)+\underline{i},
\end{equation}
which appears in $\S$ 1.2 of \cite{TS21}.
For $n\geq 1$,
define the subspace 
$V_{n}\equiv {\rm Lin}\langle\{e_{1},\ldots,e_{n}\}\rangle $
of $\ltn$.
From (\ref{eqn:newramatrix}),
$R^{(n,m)}(V_{n}\otimes V_{m})=V_{n}\otimes V_{m}$.
Hence 
$R^{(n,m)}|_{V_{n}\otimes V_{m}}$
is a linear transformation on 
$V_{n}\otimes V_{m}$. Especially,
it is a unitary operator on $V_{n}\otimes V_{m}$. 

Here we consider the case $(n,m)=(2,3)$.
We consider the restriction of the local $R$-matrix  $R^{(2,3)}$ 
on the finite dimensional subspace $V_{2}\otimes V_{3}$.
For example,
%
%
\begin{equation}
\label{eqn:rtwothree}
R^{(2,3)}(e_{1}\otimes e_{1})=e_{1}\otimes e_{1},
\quad 
R^{(2,3)}(e_{1}\otimes e_{2})=e_{2}\otimes e_{1}.
\quad 
R^{(2,3)}(e_{2}\otimes e_{1})
=e_{2}\otimes e_{2}.
\end{equation}
Define the projection
$E:\ltn\otimes \ltn\to V_{2}\otimes V_{3}$
and 
let $e_{i,j}\equiv e_{i}\otimes e_{j}$ in $V_{2}\otimes V_{3}$
for $(i,j)\in \{1,2\}\times \{1,2,3\}$.
With respect to the ordering of the basis 
$e_{1,1},e_{1,2},e_{1,3},e_{2,1},e_{2,2},e_{2,3}$
of $V_{2}\otimes V_{3}$,
the linear transformation 
$ER^{(2,3)}E$ is represented as the following $6\times 6$-matrix:
%
%
\begin{equation}
\label{eqn:matrixfour}
ER^{(2,3)}E=
\left[
\begin{array}{cccccc}
1 &&&&&\\
&&1 &&&\\
&&&&1&\\
& 1 &&&&\\
&&&1 &&
\\
&&&&& 1
\end{array}
\right].
\end{equation}
In consequence,
$ER^{(2,3)}E$ is regarded as a cyclic permutation matrix.
In this way, $R^{(2,3)}$ is not the identity operator on $\ltn\otimes \ltn$.
}
\end{ex}
%
%
\begin{ex}
\label{ex:fourth}
{\rm
We show an example which does not satisfy (\ref{eqn:implement})
because of the non-commutativity of states.
Let $\pi_{2}$ and $\omega_{2}$ 
be as in 
(\ref{eqn:pin}) and (\ref{eqn:vectorone}),
respectively.
Define $\overline{\omega}_{2}\equiv \omega_{2}\circ \alpha$
where $\alpha$ denotes the flip automorphism of the canonical generators 
$s^{(2)}_{1},s^{(2)}_{2}$ of $\co{2}$:
%
%
\begin{equation}
\label{eqn:flip}
\alpha(s_{1}^{(2)})=s_{2}^{(2)},\quad \alpha(s_{2}^{(2)})=s_{1}^{(2)},
\end{equation}
and let $\overline{\pi}_{2}$ denote the GNS representation of $\co{2}$ by 
$\overline{\omega}_{2}$ which is identified with $\pi_{2}\circ \alpha$.
Then $\omega_{2}\star \overline{\omega}_{2}\ne 
\overline{\omega}_{2}\star \omega_{2}$.
Let $v\equiv e_{1}\otimes e_{1}$.
Then $(\pi_{2}\otimes \overline{\pi}_{2})(\delp(s_{2}^{(4)}))v=v$.
Define the operator $R$ on $\ltn\otimes \ltn$ by
%
%
\begin{equation}
\label{eqn:rpi}
R(\pi_{2}\otimes \overline{\pi}_{2})(\delp(x))v
\equiv 
(\pi_{2}\otimes \overline{\pi}_{2})(\delp^{op}(x))v\quad(x\in\co{4}).
\end{equation}
Then
%
%
\begin{equation}
\label{eqn:rpitwo}
R(\pi_{2}\otimes \overline{\pi}_{2})(\delp(s_{2}^{(4)}))R^{*}v
=R(\pi_{2}\otimes \overline{\pi}_{2})(\delp(s_{2}^{(4)}))v=Rv=v.
\end{equation}
Since
$(\pi_{2}\otimes \overline{\pi}_{2})(\delp^{op}(s_{2}^{(4)}))
=\pi_{2}(s_{2}^{(2)})\otimes \overline{\pi}_{2}(s_{1}^{(2)})$
and 
$\langle v|
(\pi_{2}(s_{2}^{(2)})\otimes \overline{\pi}_{2}(s_{1}^{(2)}))v
\rangle =0$, we see that
$v\ne (\pi_{2}(s_{2}^{(2)})\otimes \overline{\pi}_{2}(s_{1}^{(2)}))v$.
Hence 
$R(\pi_{2}\otimes \overline{\pi}_{2})(\delp(s_{2}^{(4)}))R^{*}\ne 
(\pi_{2}\otimes \overline{\pi}_{2})(\delp^{op}(s_{2}^{(4)}))$.
Hence $R$ does not satisfy (\ref{eqn:implement}).
}
\end{ex}

%
%

%
\end{document}